\documentclass[10pt,letterpaper]{book}
\usepackage{amsmath,amsfonts,theorem}
\newtheorem{theorem}{Theorem}
\newtheorem{proposition}[theorem]{Proposition}
\renewcommand\chaptermark[1]{\markboth{Vector-valued quadratic forms}
{Vector-valued quadratic forms}}
\renewcommand\sectionmark[1]{\markboth{Vector-valued quadratic forms}
{Vector-valued quadratic forms}}
\begin{document}
\chapter*{
  {\vspace{-5ex}\hfill \texttt{\begin{minipage}{.7\textwidth}\small Submitted to Workshop on Open Problems in Mathematical Systems and Control Theory, MTNS'02\end{minipage}}}
  \\ [3ex]Computational problems for vector-valued quadratic forms}
\chaptermark{Vector-valued quadratic forms}
\sectionmark{Vector-valued quadratic forms}

\setcounter{chapter}{1}

\begin{minipage}[t]{0.46\hsize}\raggedright
{\Large Francesco Bullo}\\
Coordinated Science Laboratory\\
University of Illinois\\
Urbana-Champaign, IL 61801\\
United States\\
\texttt{bullo@uiuc.edu}
\end{minipage}\hspace{0.07\hsize}
\begin{minipage}[t]{0.46\hsize}\raggedright
{\Large Jorge Cort{\'e}s}\\
Systems, Signals and Control\\
University of Twente\\
Enschede, 7500 AE\\
The Netherlands\\
\texttt{j.cortesmonforte@utwente.nl}
\end{minipage}\\[1em]
\begin{minipage}[t]{0.46\hsize}\raggedright
{\Large Andrew D.~Lewis}\\
Mathematics \& Statistics\\
Queen's University\\
Kingston, ON K7L 3N6\\
Canada\\
\texttt{andrew@mast.queensu.ca}
\end{minipage}\hspace{0.07\hsize}
\begin{minipage}[t]{0.46\hsize}\raggedright
{\Large Sonia Mart{\'\i}nez}\\
Institute of Mathematics \& Physics\\
High Council of Scientific Research\\
Madrid, 28006\\
Spain\\
\texttt{s.martinez@imaff.cfmac.csic.es}
\end{minipage}

\addtocontents{toc}{\textit{Francesco Bullo, Jorge Cort{\'e}s, Andrew
D.~Lewis, Sonia Mart{\'\i}nez\\}}

\section{Problem statement and historical remarks}

For $\mathbb{R}$-vector spaces $U$ and $V$ we consider a symmetric bilinear
map $B\colon U\times U\rightarrow V$\@.  This then defines a quadratic map
$Q_B\colon U\rightarrow V$ by $Q_B(u)=B(u,u)$\@.  Corresponding to each
$\lambda\in V^*$ is a $\mathbb{R}$-valued quadratic form $\lambda Q_B$ on $U$
defined by $\lambda Q_B(u)=\lambda\cdot Q_B(u)$\@.  $B$ is \emph{definite} if
there exists $\lambda\in V^*$ so that $\lambda Q_B$ is positive-definite.
$B$ is \emph{indefinite} if for each $\lambda\in V^*$\@, $\lambda Q_B$ is
neither positive nor negative-semidefinite.  The problem we consider is as
follows.
\begin{trivlist}\item[]\it
Given a symmetric bilinear map $B\colon U\times U\rightarrow V$\@:
\begin{enumerate}
\item are there necessary and sufficient conditions, checkable in
polynomial-time, for determining when $Q_B$ is surjective?
\item if $Q_B$ is surjective, given $v\in V$ is there a polynomial-time
algorithm for finding a point $u\in Q_B^{-1}(v)$\@?
\item are there necessary and sufficient conditions, checkable in
polynomial-time, for determining when $B$ is indefinite?
\end{enumerate}
\end{trivlist}

Before we comment on how our problem impinges on control theory, let us
provide some historical context for it as a purely mathematical one.  The
classification of $\mathbb{R}$-valued quadratic forms is well understood.
However, for quadratic maps taking values in vector spaces of dimension two
or higher, the classification problem becomes more difficult.  The theory can
be thought of as beginning with the work of Kronecker, who obtained a finite
classification for pairs of symmetric matrices.  For three or more symmetric
matrices, that the classification problem has an uncountable number of
equivalence classes for a given dimension of the domain follows from the work
of Kac~\cite{VGK:83}\@.  For quadratic forms, in a series of papers Dines
(see~\cite{LLD:43} and references cited therein) investigated conditions when
a finite collection of $\mathbb{R}$-valued quadratic maps were simultaneously
positive-definite.  The study of vector-valued quadratic maps is ongoing. A
recent paper is~\cite{DBL/LMS:99}\@, to which we refer for other references.

\section{Control theoretic motivation}

Interestingly and perhaps not obviously, vector-valued quadratic forms come
up in a variety of places in control theory.  We list a few of these here.

\vspace{-0.2em}
\paragraph{Optimal control:} Agra\v{c}hev~\cite{AAA:90} explicitly
realises second-order conditions for optimality in terms of vector-valued
quadratic maps.  The geometric approach leads naturally to the consideration
of vector-valued quadratic maps, and here the necessary conditions involve
definiteness of these maps.  Agra\v{c}hev and
Gamkrelidze~\cite{AAA:88,AAA/RVG:91}\@ look at the map $\lambda\mapsto\lambda
Q_B$ from $V^*$ into the set of vector-valued quadratic maps.  Since $\lambda
Q_B$ is a $\mathbb{R}$-valued quadratic form, one can talk about its
\emph{index} and \emph{rank} (the number of $-1$'s and nonzero terms,
respectively, along the diagonal when the form is diagonalised).
In~\cite{AAA:88,AAA/RVG:91} the topology of the surfaces of constant index of
the map $\lambda\mapsto\lambda Q_B$\@ is investigated.

\vspace{-0.2em}
\paragraph{Local controllability:}  The use of vector-valued quadratic 
forms arises from the attempt to arrive at feedback-invariant conditions for
controllability. Basto-Gon{\c{c}}alves~\cite{JBG:98} gives a second-order
sufficient condition for local controllability, one of whose hypotheses is
that a certain vector-valued quadratic map be indefinite (although the
condition is not stated in this way).  This condition is somewhat refined
in~\cite{RMH/ADL:02a}\@, and a necessary condition for local controllability
is also given. Included in the hypotheses of the latter is the condition that
a certain vector-valued quadratic map be definite.

We note that Sontag~\cite{EDS:88} and Kawski~\cite{MK:90a} have shown that
the problem of determining local controllability is NP-hard.  Our problem of
asking whether there is a polynomial-time algorithm for determining the
indefiniteness of a quadratic map is not inconsistent with these results
since, in terms of controllability, our problem concerns only second-order
conditions.  However, it would be interesting if even second-order conditions
were shown to be difficult computationally.

\vspace{-0.2em}
\paragraph{Control design via power series methods and singular inversion:} 

Numerous control design problems can be tackled using power series and
inversion methods. The early references~\cite{EGA:61,AH:75} show how to solve
the optimal regulator problem and the recent work in~\cite{WTC/FB:01}
proposes local steering algorithms.  These strong results apply to linearly
controllable systems, and no general methods are yet available under only
second-order sufficient controllability conditions. While for linearly
controllable systems the classic inverse function theorem suffices, the key
requirement for second-order controllable systems is the ability to check
surjectivity and compute an inverse function for certain vector-valued
quadratic forms.

\vspace{-0.2em}
\paragraph{Dynamic feedback linearisation:} In \cite{WMS:93} Sluis gives a
necessary condition for the dynamic feedback linearisation of a system
\[
\dot{x}=f(x,u),\quad x\in\mathbb{R}^n,\ u\in\mathbb{R}^m.
\]
The condition is that for each $x\in\mathbb{R}^n$, the set $D_x=\{f(x,u)\in
T_x\mathbb{R}^n|\enspace u\in\mathbb{R}^m\}$ admits a \emph{ruling}\@, that
is, a foliation of $D_x$ by lines.  Some manipulations with differential
forms turns this necessary condition into one involving a symmetric bilinear
map $B$\@.  The condition, it turns out, is that $Q_B^{-1}(0)\not=\{0\}$\@.
This is shown by Agra\v{c}hev~\cite{AAA:88} to generically imply that $Q_B$
is surjective.

\section{Known results}

Let us state a few results along the lines of our problem statement that are
known to the authors.  The first is readily shown to be true
(see~\cite{RMH/ADL:02a} for the proof).  If $X$ is a topological space with
subsets $A\subset S\subset X$\@, we denote by $\textup{int}_S(A)$ the
interior of $A$ relative to the induced topology on $S$\@.  If $S\subset V$,
$\textup{aff}(S)$ and $\textup{conv}(S)$ denote, respectively, the affine
hull and the convex hull of $S$\@.
\begin{proposition}
Let\/ $B\colon U\times U\rightarrow V$ be a symmetric bilinear map with\/ $U$
and\/ $V$ finite-dimensional.  The following statements hold:
\begin{enumerate}
\item[(i)] $B$ is indefinite if and only if\/
$0\in\textup{int}_{\textup{aff}(\textup{image}(Q_B))}
(\textup{conv}(\textup{image}(Q_B)))$\@;
\item[(ii)] $B$ is definite if and only if there exists a hyperplane\/
$P\subset V$ so that\/ $\textup{image}(Q_B)\cap P=\{0\}$ and so that\/
$\textup{image}(Q_B)$ lies on one side of\/ $P$\@;
\item[(iii)] if\/ $Q_B$ is surjective then\/ $B$ is indefinite.
\end{enumerate}
\end{proposition}
The converse of~(iii) is false.  The quadratic map from $\mathbb{R}^3$ to
$\mathbb{R}^3$ defined by $Q_B(x,y,z)=(xy,xz,yz)$ may be shown to be
indefinite but not surjective.

Agra\v{c}hev and Sarychev~\cite{AAA/AVS:96} prove the following
result.  We denote by $\textup{ind}(Q)$ the index of a quadratic map $Q\colon
U\rightarrow\mathbb{R}$ on a vector space $U$\@.
\begin{proposition}\label{prop:sufficient}
Let\/ $B\colon U\times U\rightarrow V$ be a symmetric bilinear map with\/ $V$
finite-dimensional.  If\/ $\textup{ind}(\lambda Q_B)\ge\dim(V)$ for 
any\/ $\lambda\in V^*\setminus\{0\}$ then\/ $Q_B$ is surjective.
\end{proposition}
This sufficient condition for surjectivity is not necessary.  The quadratic
map from $\mathbb{R}^2$ to $\mathbb{R}^2$ given by $Q_B(x,y)=(x^2-y^2,xy)$ is
surjective, but does not satisfy the hypotheses of
Proposition~\ref{prop:sufficient}\@.

\section{Problem simplification}

One of the difficulties with studying vector-valued quadratic maps is that
they are somewhat difficult to get ones hands on.  However, it turns out to
be possible to simplify their study by a reduction to a rather concrete
problem.  Here we describe this process, only sketching the details of how 
to go from a given symmetric bilinear map $B\colon U\times U\rightarrow V$ 
to the reformulated end problem.  We first simplify the problem by imposing
an inner product on $U$ and choosing an orthonormal basis so that we may take
$U=\mathbb{R}^n$\@.

We let $\textup{Sym}_n(\mathbb{R})$ denote the set of symmetric $n\times n$
matrices with entries in $\mathbb{R}$\@.  On $\textup{Sym}_n(\mathbb{R})$ we
use the canonical inner product
\[
\langle\boldsymbol{A},\boldsymbol{B}\rangle=
\textup{tr}(\boldsymbol{A}\boldsymbol{B}).
\]
We consider the map
$\pi\colon\mathbb{R}^n\rightarrow\textup{Sym}_n(\mathbb{R})$ defined by
$\pi(\boldsymbol{x})=\boldsymbol{x}\boldsymbol{x}^t$\@, where $^t$ denotes
transpose.  Thus the image of $\pi$ is the set of symmetric matrices of rank
at most one.  If we identify
$\textup{Sym}_n(\mathbb{R})\simeq\mathbb{R}^n\otimes\mathbb{R}^n$\@, then
$\pi(\boldsymbol{x})=\boldsymbol{x}\otimes\boldsymbol{x}$\@.  Let $K_n$ be
the image of $\pi$ and note that it is a cone of dimension $n$ in
$\textup{Sym}_n(\mathbb{R})$ having a singularity only at its vertex at the
origin.  Furthermore, $K_n$ may be shown to be a subset of the hypercone in
$\textup{Sym}_n(\mathbb{R})$ defined by those matrices $\boldsymbol{A}$ in
$\textup{Sym}_n(\mathbb{R})$ forming angle $\arccos(\frac{1}{n})$ with the
identity matrix.  Thus the ray from the origin in
$\textup{Sym}_n(\mathbb{R})$ through the identity matrix is an axis for the
cone $K_N$\@.  In algebraic geometry, the image of $K_n$ under the
projectivisation of $\textup{Sym}_n(\mathbb{R})$ is known as the
\emph{Veronese surface}~\cite{JH:92}\@, and as such is well-studied, although
perhaps not along lines that bear directly on the problems of interest in
this article.

We now let $B\colon\mathbb{R}^n\times\mathbb{R}^n\rightarrow V$ be a
symmetric bilinear map with $V$ finite-dimensional. Using the universal
mapping property of the tensor product, $B$ induces a \emph{linear} map
$\tilde{B}\colon\textup{Sym}_n(\mathbb{R})\simeq
\mathbb{R}^n\otimes\mathbb{R}^n\rightarrow V$ with the property that
$\tilde{B}\circ\pi=B$\@.  The dual of this map gives an injective linear map
$\tilde{B}^*\colon V^* \rightarrow\textup{Sym}_n(\mathbb{R})$ (here we assume
that the image of $B$ spans $V$).  By an appropriate choice of inner product
on $V$ one can render the embedding $\tilde{B}^*$ an isometric embedding of
$V$ in $\textup{Sym}_n(\mathbb{R})$\@.  Let us denote by $L_B$ the image of
$V$ under this isometric embedding.  One may then show that with these
identifications, the image of $Q_B$ in $V$ is the orthogonal projection of
$K_n$ onto the subspace $L_B$\@.  Thus we reduce the problem to one of
orthogonal projection of a canonical object, $K_n$\@, onto a subspace in
$\textup{Sym}_n(\mathbb{R})$\@!  To simplify things further, we decompose
$L_B$ into a component along the identity matrix in
$\textup{Sym}_n(\mathbb{R})$ and a component orthogonal to the identity
matrix.  However, the matrices orthogonal to the identity are readily seen to
simply be the traceless $n\times n$ symmetric matrices.  Using our picture of
$K_n$ as a subset of a hypercone having as an axis the ray through the
identity matrix, we see that questions of surjectivity, indefiniteness, and
definiteness of $B$ impact only on the projection of $K_n$ onto that
component of $L_B$ orthogonal to the identity matrix.

The following summarises the above discussion.
\begin{trivlist}\item[]\it
The problem of studying the image of a vector-valued quadratic form can be
reduced to studying the orthogonal projection of\/
$K_n\subset\textup{Sym}_n(\mathbb{R})$\@, the unprojectivised Veronese
surface, onto a subspace of the space of traceless symmetric matrices.
\end{trivlist}
This is, we think, a beautiful interpretation of the study of vector-valued
quadratic mappings, and will surely be a useful formulation of the problem.
For example, with it one easily proves the following result.
\begin{proposition}
If\/ $\dim(U)=\dim(V)=2$ with\/ $B\colon U\times U\rightarrow V$ a symmetric
bilinear map, then\/ $Q_B$ is surjective if and only if\/ $B$ is indefinite.
\end{proposition}

\section*{References}
\begin{list}{[\arabic{enumiv}]}
  {\settowidth\labelwidth{[10]}\leftmargin\labelwidth
   \advance\leftmargin\labelsep\usecounter{enumiv}}

\bibitem{AAA:88}
A.~A.~Agra\v{c}hev.
\newblock The topology of quadratic mappings and {H}essians of smooth mappings.
\newblock {\em J. Soviet Math.}, 49(3):990--1013, 1990.

\bibitem{AAA:90}
A.~A.~Agra\v{c}hev.
\newblock Quadratic mappings in geometric control theory.
\newblock {\em J. Soviet Math.}, 51(6):2667--2734, 1990.

\bibitem{AAA/RVG:91}
A.~A.~Agra\v{c}hev and R.~V.~Gamkrelidze.
\newblock Quadratic mappings and vector functions: {E}uler characteristics of
  level sets.
\newblock {\em J. Soviet Math.}, 55(4):1892--1928, 1991.

\bibitem{AAA/AVS:96}
A.~A.~Agra\v{c}hev and A.~V.~Sarychev.
\newblock Abnormal sub-{R}iemannian geodesics: {M}orse index and rigidity.
\newblock {\em Ann. Inst. H. Poincar\'e. Anal. Non Lin\'eaire}, 13(6):635--690,
  1996.

\bibitem{EGA:61}
{\`E}.~G. Al'brekht.
\newblock On the optimal stabilization of nonlinear systems.
\newblock {\em J. Appl. Math. and Mech.}, 25:1254--1266, 1961.

\bibitem{JBG:98}
J.~{Basto-Gon{\c{c}}alves}.
\newblock Second-order conditions for local controllability.
\newblock {\em Systems Control Lett.}, 35(5):287--290, 1998.

\bibitem{WTC/FB:01}
W.~T.~Cerven and F.~Bullo.
\newblock Constructive controllability algorithms for motion planning and
  optimization.
\newblock Preprint, November 2001.

\bibitem{LLD:43}
L.~L.~Dines.
\newblock On linear combinations of quadratic forms.
\newblock {\em Bull. Amer. Math. Soc. (N.S.)}, 49:388--393, 1943.

\bibitem{AH:75}
A.~Halme.
\newblock On the nonlinear regulator problem.
\newblock {\em J. Optim. Theory Appl.}, 16(3-4):255--275, 1975.

\bibitem{JH:92}
J.~Harris.
\newblock {\em Algebraic Geometry: A First Course}.
\newblock Number 133 in Graduate Texts in Mathematics. Springer-Verlag, New
  York-Heidelberg-Berlin, 1992.

\bibitem{RMH/ADL:02a}
R.~M.~Hirschorn and A.~D.~Lewis.
\newblock Second-order controllability conditions with weak hypotheses for
  control affine systems.
\newblock Preprint, February 2002.

\bibitem{VGK:83}
V.~G.~Kac.
\newblock Root systems, representations of quivers and invariant theory.
\newblock In {\em Invariant theory}, number 996 in Lecture Notes in
  Mathematics, pages 74--108. Springer-Verlag, New York-Heidelberg-Berlin,
  1983.

\bibitem{MK:90a}
M.~Kawski.
\newblock The complexity of deciding controllability.
\newblock {\em Systems Control Lett.}, 15(1):9--14, 1990.

\bibitem{DBL/LMS:99}
D.~B.~Leep and L.~M.~Schueller.
\newblock Classification of pairs of symmetric and alternating bilinear forms.
\newblock {\em Exposition. Math.}, 17(5):385--414, 1999.

\bibitem{WMS:93}
W.~M.~Sluis.
\newblock A necessary condition for dynamic feedback linearization.
\newblock {\em Systems Control Lett.}, 21(4):277--283, 1993.

\bibitem{EDS:88}
E.~D.~Sontag.
\newblock Controllability is harder to decide than accessibility.
\newblock {\em SIAM J. Control Optim.}, 26(5):1106--1118, 1988.

\end{list}

\end{document}